\documentclass{amsart}

\usepackage{graphicx}
\usepackage{amsmath}%
\usepackage{amsfonts}%
\usepackage{amssymb}

\usepackage{epsfig}

\usepackage{psfrag}

\newtheorem{theorem}{Theorem}[section]

\newtheorem{corollary}[theorem]{Corollary}

\newtheorem{lemma}[theorem]{Lemma}

\newtheorem{proposition}[theorem]{Proposition}

\numberwithin{equation}{section}
\newcommand{\RR}{{\mathbb{R}}}

\def\ova{{\bf a}}
\def\ba{\ova}
\def\lam{\lambda}
\def\Th{\theta}
\def\gam{\gamma}
\def\Om{\Omega}

\def\Gam{\Gamma}

\def\sig{\sigma}
\def\Uk{{\mathcal U}}
\def\ela{\stackrel{\lam}{=}}

\def\supp{{\rm supp}}

\def\nuhat{\widehat{\nu}}
\def\muhat{\widehat{\mu}}

\newcommand{\R}{{\mathbb R}} \newcommand{\E}{{\mathbb E}\,} \newcommand{\Prob}{{\mathbb P}\,} \newcommand{\Nat}{{\mathbb N}}  \newcommand{\Ak}{{\mathcal A}}
\def\Pk{{\mathcal P}}
\def\Gk{{\mathcal G}}
\def\Z{{\mathbb Z}}
\newcommand{\Ek}{{\mathcal E}} \newcommand{\Fk}{{\mathcal F}} \def\Mk{{\mathcal M}}  \newcommand{\eps}{{\varepsilon}} \newcommand{\es}{\emptyset} \newcommand{\half}{{\frac{1}{2}}}

\begin{document}

\title[BRW with Exponentially Decreasing Steps]{Branching Random Walk with Exponentially Decreasing Steps, and Stochastically Self-Similar Measures}

\author{Itai Benjamini \and Ori Gurel-Gurevich \and Boris Solomyak}

\address{Boris Solomyak, Box 354350, Department of Mathematics,
University of Washington, Seattle WA 98195, USA}
\email{solomyak@math.washington.edu}
\address{Itai Benjamini and Ori Gurel-Gurevich,
Department of Theoretical Mathematics, Weizmann Institute of Science,
Rehovot, 76100, Israel}

\thanks{
\\ \indent
Research of Solomyak was partially supported by
NSF grant DMS 0355187.
}

\begin{abstract}
We consider a Branching Random Walk on $\RR$ whose step size decreases by a
fixed factor, $0<\lambda<1$, with each turn. This process generates
a random probability measure on $\RR$, that is, the limit of uniform
distribution among the $2^n$ particles of the $n$-th step.
We present an initial investigation of the limit measure and its support.
We show, in particular, that (1) for almost every $\lam>1/2$ the limit
measure is almost surely (a.s.) absolutely continuous with respect to
the Lebesgue
measure, but for Pisot $1/\lam$ it is a.s.\ singular;
(2) for all $\lam> (\sqrt{5}-1)/2$ the support of the measure is a.s.\ the closure of its interior;
(3) for Pisot $1/\lam$ the support of the measure is ``fractured'':
it is a.s.\ disconnected and the components of the complement are not isolated on both sides.
\end{abstract}

\maketitle

\section{Introduction and statement of results} 
A \emph{Branching Random Walk (BRW)} on $X$ is a
random map from the complete infinite binary tree, $T=\{1,2\}^*$
into $X$. We shall consider a symmetric BRW on $\RR$ with
exponentially decreasing steps, defined as follows. Start with a
single particle at 0. At each step each particle multiplies to two
particles, and each independently takes a step of size $\lambda^{n-1}$
to either direction with equal probabilities.  Another equivalent
useful formulation is this: for each vertex $v$ of the binary tree
let $a_v$ be equal to $+1$ or $-1$ with probability $\half$
independently. Using these lotteries we define the BRW function to
be $f(v)=\sum_{n=1}^{|v|} a_{v|n} \lambda^n$ where $v|n=v_1 \ldots v_n$.
This function is
extended to the boundary of the tree, $\partial T =\{1,2\}^\Nat$, in
the obvious manner, $f(v)=\sum_{n=1}^\infty a_{v|n} \lambda^n$. Let
$m=(\half,\half)^\Nat$ 
be the standard uniform measure on $\partial T$. Define $\mu$ to
be the image of $m$ under $f$, i.e. $\mu(E)=m(f^{-1}(E))$. This is
a random measure on the line, which depends on the choice of signs
on the tree. We are also interested in the properties of 
the compact support $S$ of $\mu$ which is clearly the image of $f$:
$S = f(\partial T)$.

One can view $\mu$ as a 
stochastically self-similar measure and $S$ as a stochastically self-similar
set with respect to the appropriate transformations.
More precisely, let $F_1(x) = \lam x+\lam$ and $F_2(x) = \lam x-\lam$, and
let 
$$
\Fk:= \frac{1}{4} \Bigl(\delta_{(F_1,F_1)} + \delta_{(F_1,F_2)} + 
\delta_{(F_2,F_1)}
+ \delta_{(F_2,F_2)}\Bigr)
$$
be a distribution on the pairs of similitudes in $\R$. Then
$$
\mu \stackrel{d}{=} \sum_{i=1}^2 \frac{1}{2} \left( \mu^{(i)} \circ
\Fk_i^{-1}\right),
$$
where $(\Fk_1,\Fk_2)$ is a random vector of similitudes distributed
according to $\Fk$, and $\mu^{(i)},\ i=1,2$, are i.i.d.\ copies of
$\mu$ independent of $(\Fk_1,\Fk_2)$. The symbol $\stackrel{d}{=}$
denotes equality in distribution.

\smallskip

There is a large literature on 
stochastically self-similar sets and measures: 
Falconer \cite{Falc}, Graf \cite{Graf},
Mauldin and Williams \cite{MW} investigated random fractal sets, and
U. Z\"ahle \cite{Zae}, Patzshke and U. Z\"ahle \cite{PZ}, Arbeiter 
\cite{Arb1,Arb2}, Olsen \cite{Olsen}, and
Hutchinson and R\"uschendorf \cite{HR} developed the theory of
random fractal measures. 
Existence, uniqueness, and convergence results have been established under
very general assumptions, but results on dimension were obtained mostly under
some separation (``non-overlapping'') conditions.
Our case is inherently overlapping for every
$\lam>0$, since
there is a positive probability 
of having $\Fk_1 = \Fk_2$.
Overlapping is also allowed in \cite[Prop.\,6.4]{Arb1}, 
where 
the translation parts of the similitudes have i.i.d.\ 
absolutely continuous distributions with
a bounded density, and dimension formulas are obtained which hold a.s.
In the recent work by T. Jordan, M. Pollicott, and K. Simon \cite{JPS},
stochastically self-affine sets and measures (with overlaps)
are studied, which of course,
includes stochastically self-similar ones as a special case. A.s.\ formulas
for the dimension and a.s.\ absolute continuity are established there under
appropriate assumptions. Both in \cite{JPS} and in the earlier work
\cite{PSS}, where a different ``overlapping'' random model was investigated,
the distributions of the vectors of similitudes are absolutely continuous.
In our model, on the other hand, these distributions are discrete.
This puts it closer to the infinite Bernoulli convolution measures,
extensively studied since the 1930's (see \cite{sixty}).

Of course, our model can be generalized in many ways: instead of the binary
tree one can consider an arbitrary rooted tree, instead of the random variables
with values $\pm 1$ one can take more general discrete random variables,
and one can consider projections of different measures on the
boundary of the tree. Other possible generalizations are mentioned in Section 5.

In Section 2 we adopt an ``intermediate-general'' viewpoint:
$T=\{1,\ldots,\ell\}^*$ is the $\ell$-regular tree for $\ell\ge 2$.
(On a deterministic non-regular tree we loose
stochastic self-similarity; although some results extend to that case,
we don't consider it here.)
The random variables $\{a_v\}_{v\in T}$
at the vertices are i.i.d.\ 
with a discrete distribution 
$\eta= \sum_{d\in D} p_d \delta_d$. Here $D\subset \R$ is a
finite set, which will be called a set of ``digits.''
First we show the ``pure types law'': for a fixed $\lam$,
the measure $\mu$ is either absolutely continuous (a.c.) with respect to the
Lebesgue measure, or purely singular, almost surely.
This is a simple consequence of uniqueness.
Then we adapt the approach of Bluhm \cite{Bluhm} to obtain estimates of
the expectation of some quantities which involve $|\muhat(t)|^2$ 
in terms of the corresponding quantities for the associated
deterministic self-similar measure. More precisely, consider the probability
measure $\nu$, which is the unique solution of the equation
\begin{equation} \label{eq-ssm}
\nu = \sum_{d\in D} p_d (\nu \circ F_d^{-1}),
\ \ \mbox{where}\ F_d(x) = \lam (x + d),
\end{equation}
see \cite{Hutch}.
It is easy to see that $\nu$ is the distribution of the random sum
$\sum_{n=1}^\infty b_n \lam^{n}$ where the coefficients $b_n$ are
i.i.d.\ with the distribution $\eta$. It follows that the Fourier transform
of $\nu$ may be computed as follows:
\begin{equation} \label{eq-fourier}
\widehat{\nu}(t) = \prod_{n=1}^\infty \widehat{\eta}(t\lam^{n}).
\end{equation}
The classical Bernoulli convolution arises this way if we take $D=\{0,1\}$
and $p_0=p_1=\half$ (or any other two distinct digits).

\begin{theorem} \label{thm-fourier} Let $\mu$ be the BRW with steps of
size $\lam^n$, with the i.i.d.\ random variables on the $\ell$-regular tree
distributed as $\eta$. Let $\nu$ be the deterministic self-similar measure
given by (\ref{eq-ssm}). Then

{\bf (i)} 
$\E|\widehat{\mu}(t)|^2 \ge |\widehat{\nu}(t)|^2$ for all $t\in \R$;

{\bf (ii)} For any $\gam\ge 0$, with $\lam^{1+2\gam} > \frac{1}{\ell}$, there
exist constants $C_1$ and $C_2$ such that
\begin{equation} \label{ineq2}
\E\Bigl(\int_{\R}  |\widehat{\mu}(t)|^2 |t|^{2\gam}\,dt\Bigr) \le C_1 + C_2
\int_{\R}  |\widehat{\nu}(t)|^2 |t|^{2\gam}\,dt.
\end{equation}
\end{theorem}

Following \cite{PS}, we use homogeneous Sobolev norms
$$
\|\nu\|^2_{2,\gam} = \int_{\R} |\nuhat(t)|^2 |t|^{2\gam}\,dt.
$$
Finiteness of $\|\nu\|_{2,\gam}$ for $\gam>0$ means that $\nu$ has
$\gam$ (fractional) derivatives in $L^2$.

\begin{corollary} \label{cor2} Let $\mu$ and $\nu$ be as in
Theorem~\ref{thm-fourier}.

{\bf (i)} If $\lim_{|t|\to\infty} |\widehat{\nu}(t)|\ne 0$, then
$\mu$ is a.s.\ singular.

{\bf (ii)} If $\E\int_{\R} |\muhat(t)|^2\, dt< \infty$, then
$\nu$ is absolutely continuous with a density in $L^2(\R)$.

{\bf (iii)} If $\nu$ is absolutely continuous
with a density in $L^2(\R)$ and $\lam >
\frac{1}{\ell}$, then
$\mu$ is a.s.\ absolutely continuous with a density in $L^2(\R)$.

{\bf (iv)} If $\|\nu\|_{2,\gam}<\infty$ and $\lam^{1+2\gam} > \frac{1}{\ell}$,
then $\|\mu\|_{2,\gam}<\infty$ a.s.
\end{corollary}

Using the results available for deterministic self-similar measures
and Bernoulli convolutions with overlaps (see \cite{Erd1,Erd2,Sol,sixty,ST} 
and references therein), we obtain a lot of information on the random measure
$\mu$. In particular, we have the following

\begin{corollary} 
\label{cor-meas} Suppose that $T$ is $\ell$-regular
for $\ell\ge 2$ and the random variables at the vertices are i.i.d.\
with the distribution $\eta = \frac{1}{m}(\delta_{d_1} + \cdots + \delta_{d_m})$
for $m\ge 2$,
where the digits $d_j$ have uniform spacing: $d_j = d_1 + (j-1)h$ for some
$h>0$. Let $\mu$ be the corresponding BRW with steps of size $\lam^{n}$.
Then

{\bf (i)} $\mu$ is a.s.\ singular
for $\lam = 1/\Th$, where $\Th<m$ is a Pisot number.

{\bf (ii)} $\mu$ is a.s.\ absolutely continuous
for a.e.\ $\lam \in (\max\{\frac{1}{\ell},\frac{1}{m}\},1)$;
\end{corollary}

Recall that a Pisot number is an algebraic integer $\Th>1$
whose conjugates (i.e.\ other zeros of the minimal polynomial) are
strictly less than one in absolute value.

\smallskip

\noindent {\bf Remarks.} 1.
Looking at the $n$-th level, we see that $S$ is covered by $\min\{\ell^n,m^n\}$
intervals of size $\sim \lam^n$. Thus, if 
$\lam < \max\{\frac{1}{\ell},\frac{1}{m}\}$, then
$S$ has Hausdorff dimension less than
one and hence
$\mu$ is singular (surely, not just almost surely).

2. Corollary~\ref{cor2}
shows that the a.s.\ properties of the
stochastically self-similar measure $\mu$ and those of its deterministic 
counterpart are closely related. There is a heuristic principle that
putting more randomness into the model increases the likelihood of 
absolute continuity. In our model the randomness is fairly ``mild,''
so that the number-theoretic phenomena associated with Pisot numbers
are preserved (unlike the models in \cite{PSS,JPS}).

3. Corollary~\ref{cor2}(ii) opens a possibility of applications
in the other direction; although this may be far-fetched, any progress
in the problem of determining precisely
for which $\lam$ the Bernoulli convolution
measure is absolutely continuous (see \cite{sixty}), would be interesting.

\smallskip

Finally, we state a corollary which gives additional information
for the most basic case $m=2$ and $\ell=2$, using the results
available for classical Bernoulli convolutions.
Motivated by \cite{Garsia}, we say that $\Th>1$ is a {\em Garsia number} if it
is an algebraic integer whose
minimal polynomial has all zeros
greater than one in absolute value and the constant term
$\pm 2$. Examples of such polynomials include $x^n-2$ for $n\ge 2$,
$x^{n+p} - x^n -2$ for $p,n\ge 1$ and $\max\{p,n\}\ge 2$, $x^3-2x-2$, etc.,
see \cite{Garsia}.

We write $\dim$ to denote the Hausdorff dimension.

\begin{corollary} \label{cor3}
 Let $\mu$ be the BRW with steps
of size $\lam^{n}$, with the i.i.d.\ random variables on the binary tree
distributed as $\half(\delta_0 + \delta_1)$.

{\bf (i)} Suppose $\Th=1/\lam\in (1,2)$ is a Garsia number.
Then $\mu$ is a.c.\ with a density in $L^2$ almost surely.

{\bf (ii)} There exist $a_k< 1,\ a_k\to 1$,
such that
$\mu$ is a.c.\ with a $k$ times differentiable density almost surely, for a.e.\
$\lam\in (a_k,1)$.

{\bf (iii)} There exists $C>0$ such that for all $\eps>0$ we have
$\dim\{\lam\in (\half+\eps,1):\,\mu \ \mbox{\rm is a.s.\ singular} \ \}
 < 1-C\eps$. 
\end{corollary}

It should be noted that the estimate in (iii) is the best known, but probably
not the best possible. The reciprocals of
Pisot numbers are the only known parameters in $(\half,1)$
for which the Bernoulli convolution measure is singular.

\smallskip

In Sections 3 and 4
we study the topological properties of the support, restricting
ourselves to the case of the binary tree and uniform Bernoulli (2-digit)
random variables, as in the beginning of the Introduction.
Here the case of classical Bernoulli convolutions cannot serve as a guide,
since for them the support of the measure is an interval whenever $\lam\ge 
\half$.

It is more convenient
to use the digits $0,1$ rather than $\pm 1$ (this is obtained by
 a linear change of variables). Then
all elements of $S = \supp(\mu)$ are of the form $\sum_{n=1}^\infty a_n 
\lam^{n}$ for $a_n \in \{0,1\}$.  Let $I :=
[0,\frac{\lam}{1-\lam}]$ and note that $S\subset I$. Also note that
if all infinite words in $\{0,1\}^\Nat$ can be ``read off'' the tree
$T$ from the root, then $S=I$. This is due to $\lam> \half$ and the
fact that every $x\in I$ has an expansion in base $\lam$ with digits
0,1. The expansion is, in general, non-unique, so the condition for
$S=I$ is only sufficient. 

For $x\in I$ consider the set of all infinite words giving an expansion of
$x$ in base $\lam$:
$$\Ek_\lam(x) = \Bigl\{a\in \{0,1\}^\Nat:\ x = \sum_{n=1}^\infty
a_n \lam^{n}\Bigr\}.
$$ 
Questions about the size of $\Ek_\lam(x)$ have been
studied, see Erd\H{o}s, Jo\'o and Komornik \cite{EJK},
Glendinning and Sidorov \cite{GS}, and references therein.
In particular, in \cite{EJK} it is proved
that for all $\lam > g:= \frac{\sqrt{5}-1}{2}$ and
$x\in (0,\frac{\lam}{1-\lam})$ the
set $\Ek_\lam(x)$ has the cardinality of continuum, and its Hausdorff
dimension in the natural metric on $\partial T$ is positive.
On the other hand, for $\lam<g$ there are $x\in (0,\frac{\lam}{1-\lam})$ 
having a unique expansion. More precisely, let
$\Psi_\lam:=\{x\in (0,\frac{\lam}{1-\lam}):\ \#\Ek_\lam(x)=1\}$.
In \cite{GS} it is proved, in particular, that $\Psi_\lam$
is countably infinite for $\lam\in (\beta,g)$ and 
is uncountable for $\lam \in (1/2,\beta]$, where $\beta\approx 0.559852...$
is the ``Komornik-Loreti constant'' \cite{KL}.

Returning to our problem, we note that $x\in S$ if and only if there exists
$a \in \Ek_\lam(x)$
which can be ``read off'' 
from the root of $T$.
The questions about ``hitting'' a given subset of the sequence space
by infinite words seen along the paths of a tree for a tree-indexed process
were considered by many authors. A set is called {\em nonpolar} or {\em polar}
according to whether or not it is hit with positive probability.
Thus, $x\in S$ with positive probability if and only if $\Ek_\lam(x)$ is 
nonpolar,
and by a result of Evans \cite[Th.\,2]{Evans},
this is equivalent to $\Ek_\lam(x)$ having positive logarithmic capacity in
the standard metric on $\{0,1\}^\Nat$. In particular, singletons are polar.
(This is also easy to see directly, since for a given word
$a \in \{0,1\}^\Nat$, the random subtree  consisting
of the edges along which we see the beginning of $a$, starting from the root,
is a Galton-Watson process with offspring distribution
$\frac{1}{4}\delta_0 + \half \delta_1 + \frac{1}{4}\delta_2$.
It is a critical branching process which dies out a.s.)
It follows that for
any $x\in \Psi_\lam$, and more generally, for any $x$ having at most
countable many
 expansions in base $\lam$ with digits $0$ and $1$,
almost surely $x\not\in S$.
We summarize this discussion in the following proposition.

\begin{proposition} \label{prop-discon}
Let $\lam \in (\half,g]$.
Let $\mu$ be the BRW with steps
of size $\lam^{n}$, with the i.i.d.\ random variables on the binary tree
distributed as $\half(\delta_0 + \delta_1)$.
Let $I = [0,\frac{\lam}{1-\lam}]$ and suppose that the set of
$x\in I$ for which $\Ek_\lam(x)$ has zero logarithmic capacity (e.g.\ if it is
countable) is dense in $I$. Then $S$ is totally disconnected.
\end{proposition}

\begin{corollary} \label{cor-gold}
If $\lam=g$, then $S$ is totally disconnected.
\end{corollary}

On the other hand, if $\Ek_\lam(x)$ has positive dimension, then
$x\in S$ with positive probability. Using the methods of 
Benjamini and Kesten \cite{BK}, who investigated when {\em all} words can be
seen from finitely many vertices, we obtain the following (see Section 3):
 
\begin{theorem} \label{th-int} Suppose that $\lam \in (g,1)$, where
$g = \frac{\sqrt{5}-1}{2}$ is the golden ratio. Let $\mu$ be the BRW with steps
of size $\lam^{n}$, with the i.i.d.\ random variables on the binary tree
distributed as $\half(\delta_0 + \delta_1)$.
Then $S=\supp(\mu)$ 
has nonempty interior, and is the closure of its interior, a.s.
\end{theorem}

The following easy statement is included for completeness.
Below $|S|$ denotes the Lebesgue measure of $S$.

\begin{proposition} \label{prop-easy}
Let $\mu$ be the BRW with steps
of size $\lam^n$, with the i.i.d.\ random variables on the binary tree
distributed as $\half(\delta_0 + \delta_1)$.
If $\lam \le \half$, then $|S|=0$ a.s.
\end{proposition}

It is also standard that for $\lam=\half$ the set $S$ has Hausdorff
dimension equal to one. This can be deduced directly or from 
\cite{Evans} (which implies that 
$S$ hits any subset of $I$ of positive dimension
with positive probability), combined with \cite{Lyons}, or, alternatively,
from \cite{Hawkes}.  

\smallskip

\noindent
{\bf Open questions.} For which $\lam\in (\half,g)$ does the set $S$ contain
intervals almost surely? For which $\lam\in (\half,g)$  is $S$ totally\
disconnected a.s.? 
What are the a.s.\ values of Hausdorff dimension and Lebesgue measure of $S$ for
the golden ratio case $\lam = g$?

\smallskip

We remark that there are many open problems
concerning the interior of self-similar sets, both
deterministic (see \cite{PeSo}) and random (see \cite{PSS}).

\smallskip

In Section 4 we turn to the question of connectedness of the support $S$.
It is obvious that $S$ has a positive probability of being disconnected
for any tree and any distribution $\eta$
at the vertices. But is $S$ disconnected almost surely? For some trees
this is not the case. For instance, if $T$ is a $4$-regular tree and
$\eta=\half(\delta_0 + \delta_1)$, then
there is a positive probability of seeing all the words in $\{0,1\}^\Nat$
from the root \cite[Cor.\,6.3]{BK}, hence a positive probability of having
$S = [0,\frac{\lam}{1-\lam}]$ for all $\lambda$.
In the case of the binary tree it seems that
$S$ should be a.s.\ disconnected, but we could only prove it in a special case.

\begin{theorem} \label{th-fract}
 Let $\mu$ be the BRW with steps
of size $\lam^{n}$, with the i.i.d.\ random variables on the binary tree
distributed as $\half(\delta_0 + \delta_1)$.
If $\lambda=1/\theta$, where $\theta$ is a Pisot number, then $S=\supp(\mu)$
is disconnected and
has no isolated gaps, almost surely. 
More precisely, let $(\alpha,\beta)$ be a component of $\R\setminus S$. 
If $\alpha >
-\infty$, then $(\alpha-\eps,\alpha)\not\subset S$ for every
$\eps>0$, and if $\beta< \infty$, then
$(\beta,\beta+\eps)\not\subset S$ for every $\eps>0$.
\end{theorem}

We remark that there are infinitely many $\lam$'s for which both
Theorem~\ref{th-int} and Theorem~\ref{th-fract} apply, see \cite{Pisot}.
These are $\lam \in (g,1)$ for which $1/\lam$ is Pisot. For the corresponding
random measure $\mu$ we know that (i) $\mu$ is singular a.s.; (ii) the
support of $\mu$ is a.s.\ the closure of a countable union of intervals, with 
the property that every gap is accumulated by gaps on both sides. Such
sets were called {\em M-cantorvals} by Mendes and Oliveira \cite{MO}; they
often appear as arithmetic sums of (deterministic) Cantor sets.

\smallskip

\noindent
{\bf Open questions.} For which $\lam > 1/2$ is the set $S$ a.s.\ disconnected?
If $S$ is disconnected with infinitely many gaps, 
what can be said about the distribution of
lengths of the complementary intervals?

\smallskip

\noindent
{\bf Figure.}
We show the approximate ``densities'' of the
random measure $\mu$ for different simulations and for several interesting
values of $\lam$. For comparison, we include pictures of the
corresponding (approximations to) Bernoulli convolutions.
To be precise, we generated $2^{20}$ points using a simulation of the
BRW, subdivided the interval $[0,\frac{\lam}{1-\lam}]$ into
$2^{10}=1024$ intervals of equal size, and plotted the corresponding histogram
(this explains the numbers on the axes).
The parameters $\lam =.565198\ldots$ and $\lam=2^{-1/2}$ are reciprocals
of  Garsia numbers, 
so by Corollary~\ref{cor3}(i), the random measure $\mu$ is a.c.\ with a density
in $L^2$ almost surely. The parameter $\lam=.618034\ldots$ is the reciprocal of
the golden ratio, a Pisot number, and $\lam = .754877\ldots$ is 
the reciprocal of the smallest Pisot number. Thus, the corresponding random
measures
are singular and have fractured support, even totally disconnected in the case
of the golden ratio, almost surely. (These phenomena, however, are not
visible in our figures, since the gaps are likely to be very small for 
larger values of $\lam$.) 

We emphasize again that there is a positive probability of having
disconnected support of the random measure $\mu$
for all $\lam < 1$, and this is a major difference with
Bernoulli convolutions. 

\begin{figure}[ht]
\begin{center}
$\begin{array}{ccc}
\multicolumn{1}{l}{} &
        \multicolumn{1}{l}{} &
        \multicolumn{1}{l}{} \\
\epsfxsize=1.7in
\epsffile{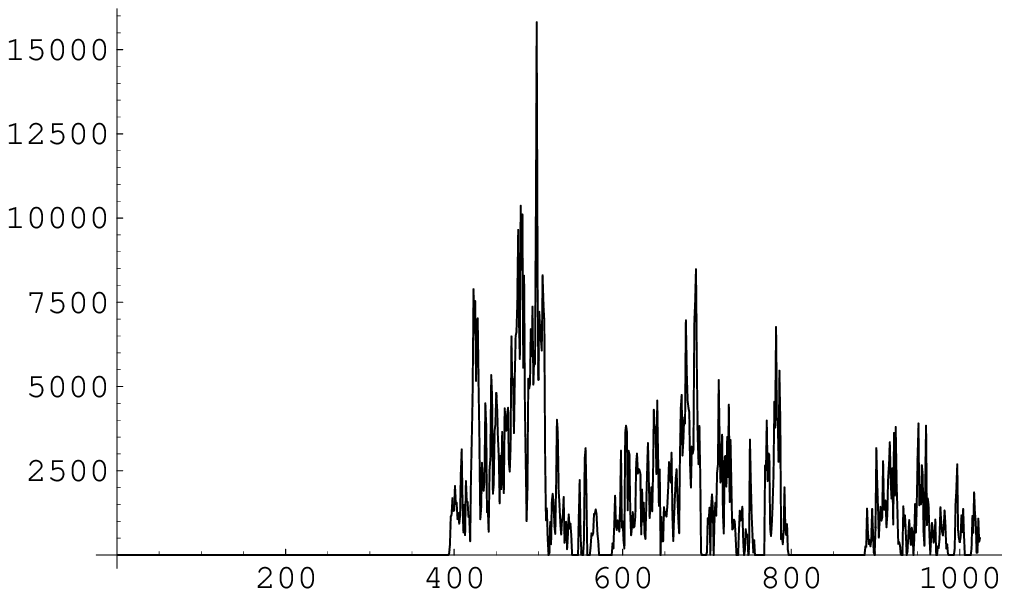} &
        \epsfxsize=1.7in
        \epsffile{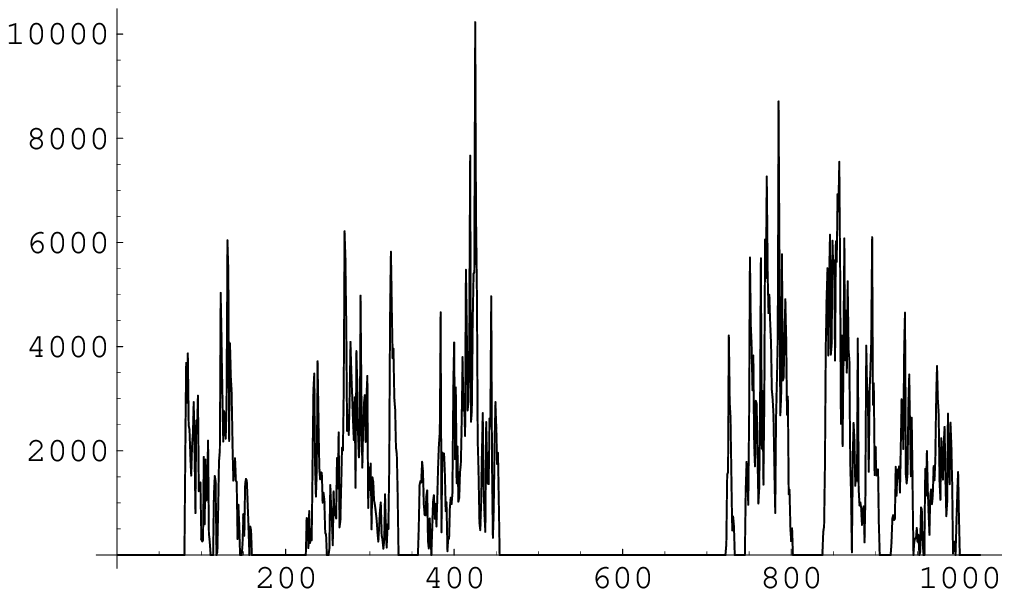}  &
        \epsfxsize=1.7in
        \epsffile{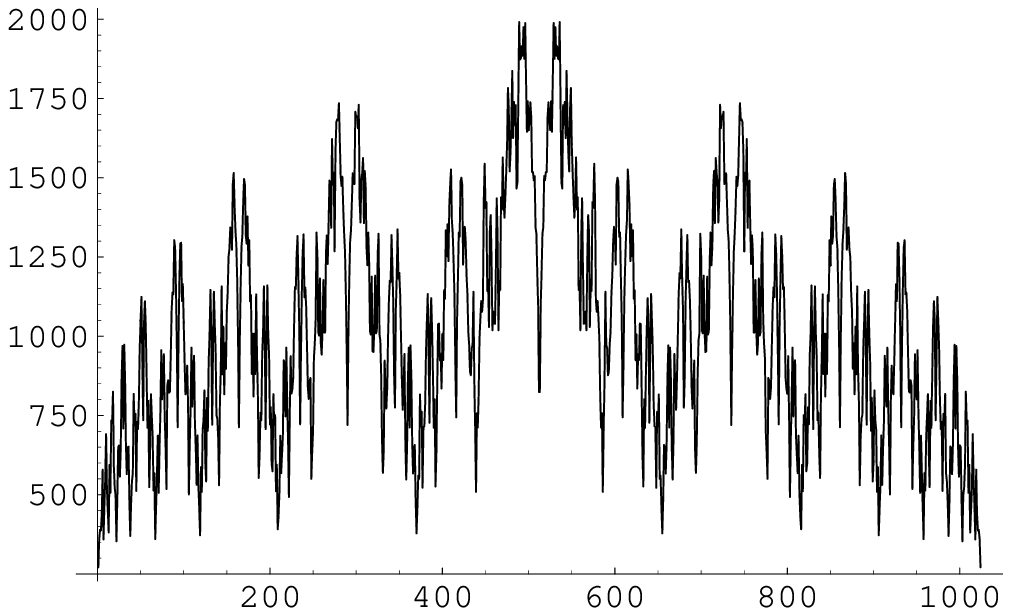} \\ & \lam=.565198\ldots
 & \\[.5cm]
\epsfxsize=1.7in
\epsffile{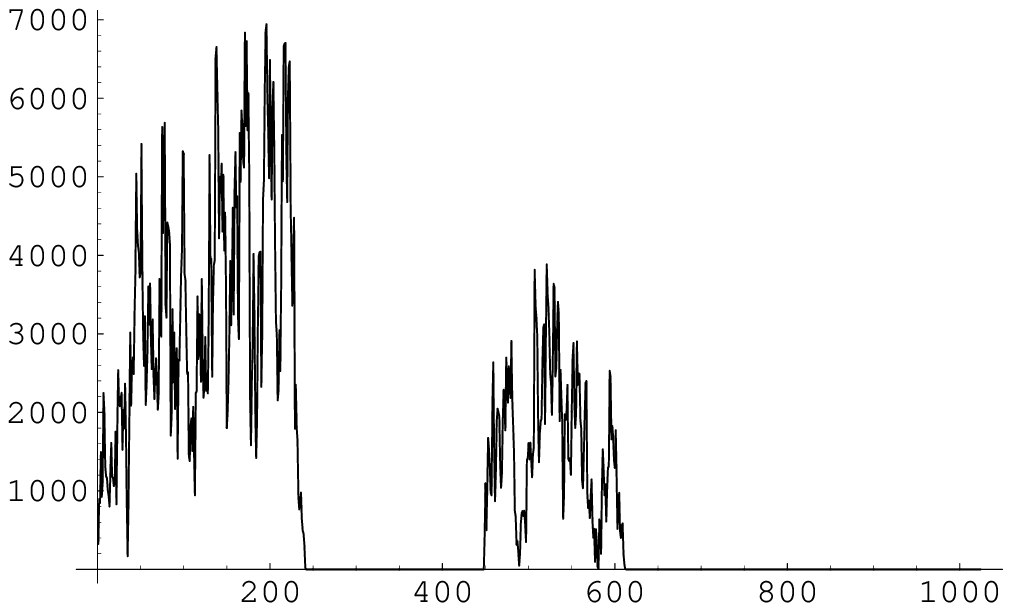} &
        \epsfxsize=1.7in
        \epsffile{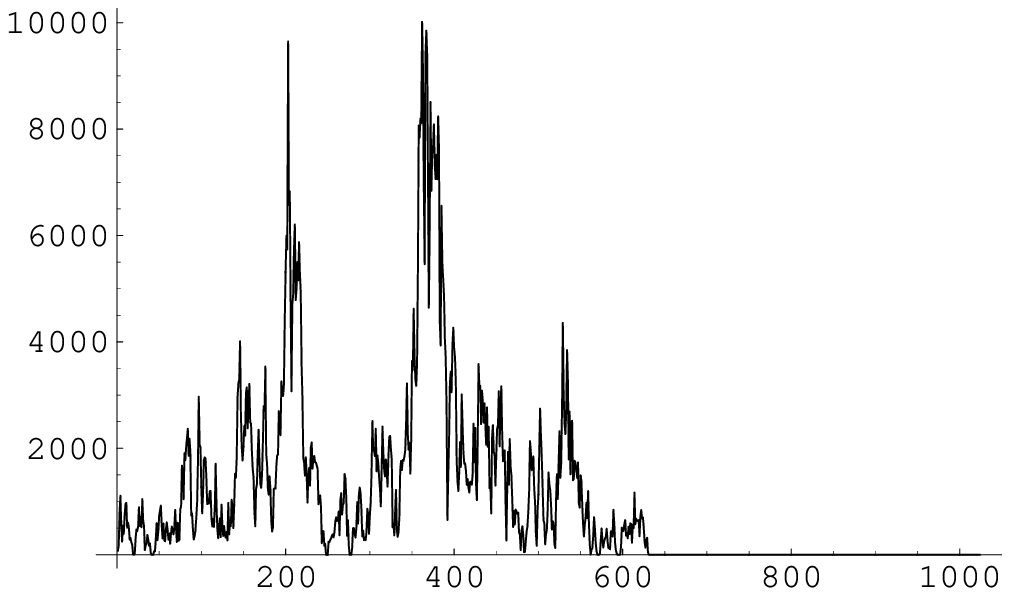}  &
        \epsfxsize=1.7in
        \epsffile{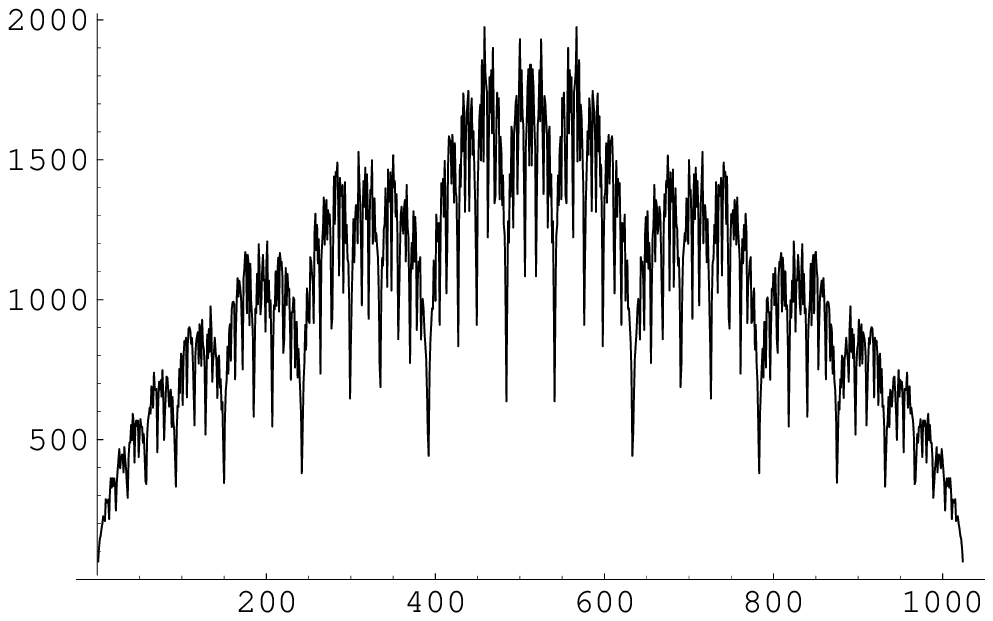}
         \\ &  \lam=.618034\ldots &  \\ [.5cm]
\epsfxsize=1.7in
\epsffile{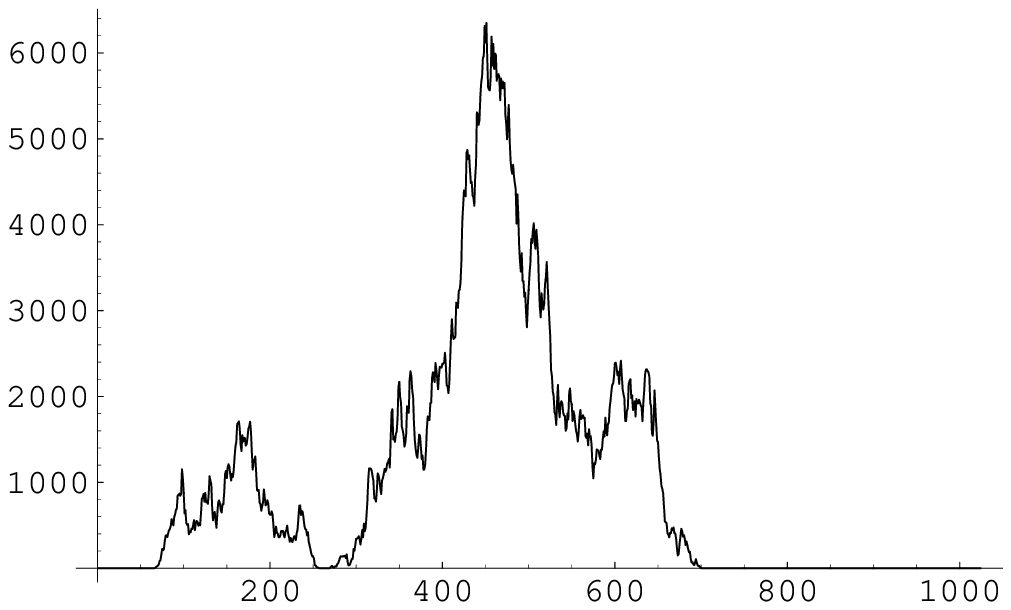} &
        \epsfxsize=1.7in
        \epsffile{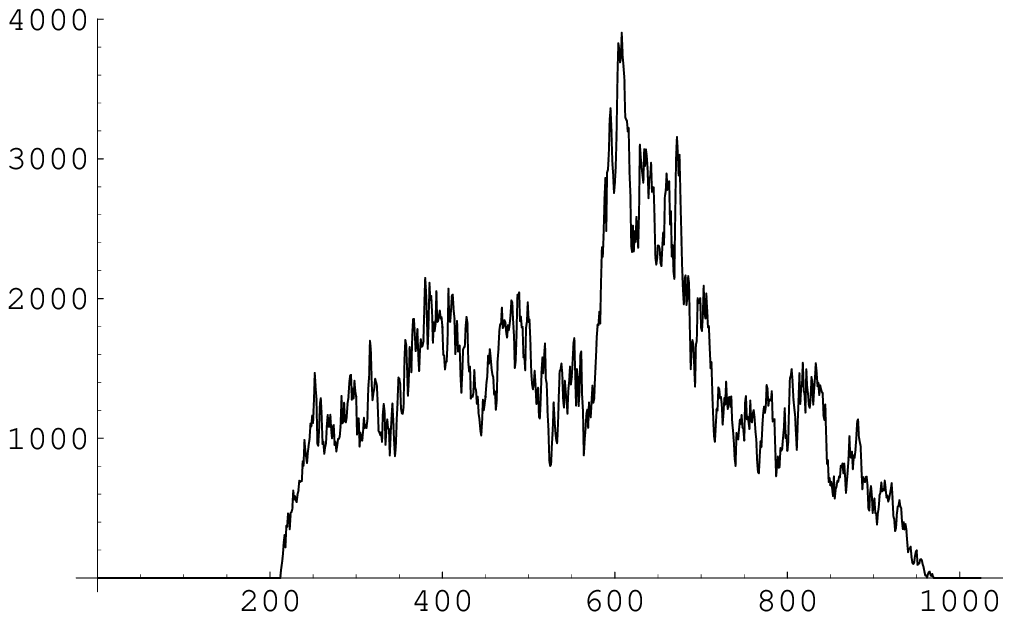}  &
        \epsfxsize=1.7in
        \epsffile{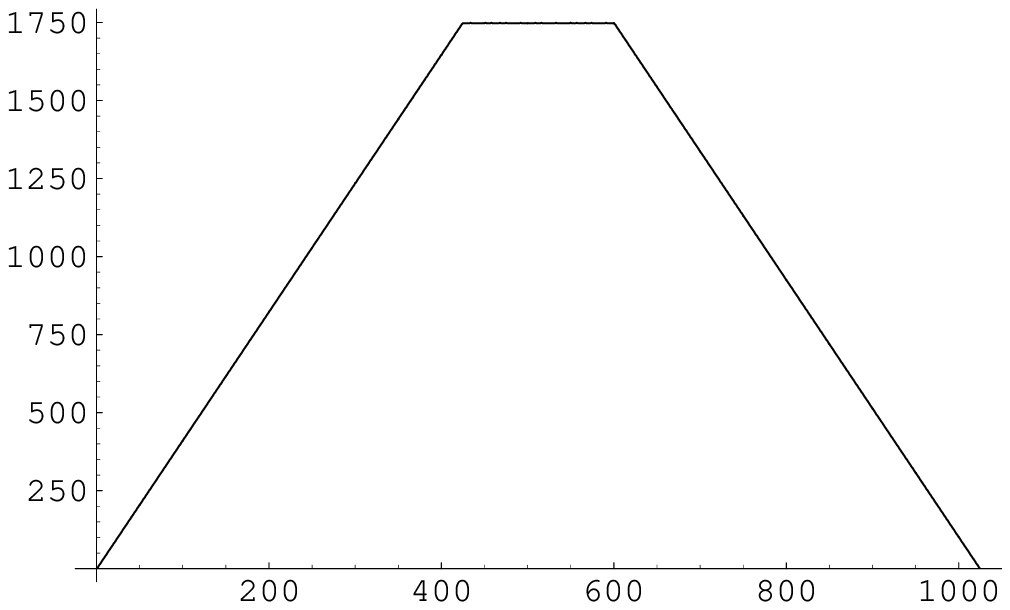}
         \\ &  \lam=2^{-1/2} &  \\[.5cm]
\epsfxsize=1.7in
\epsffile{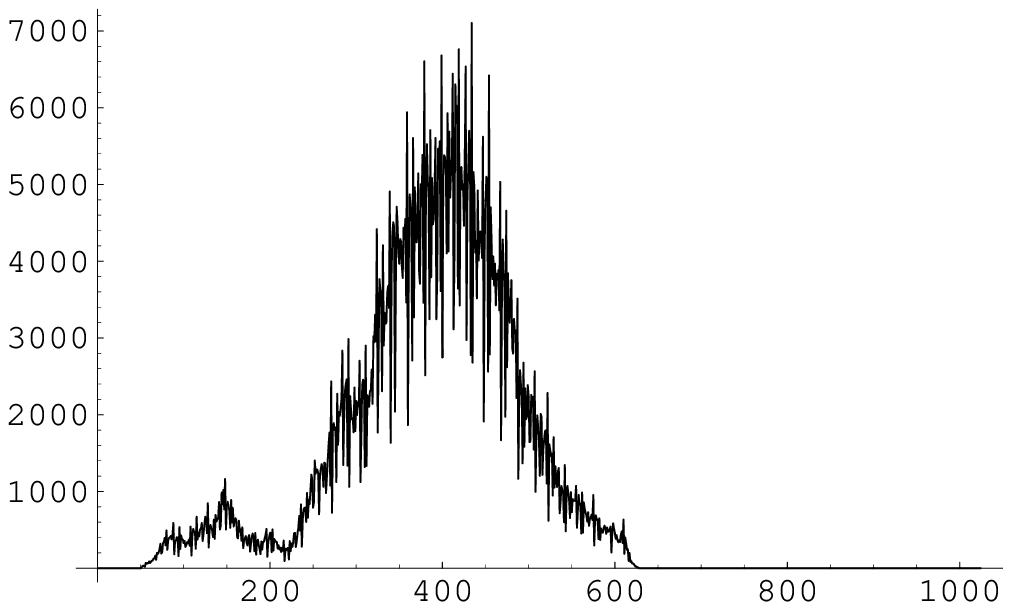} &
        \epsfxsize=1.7in
        \epsffile{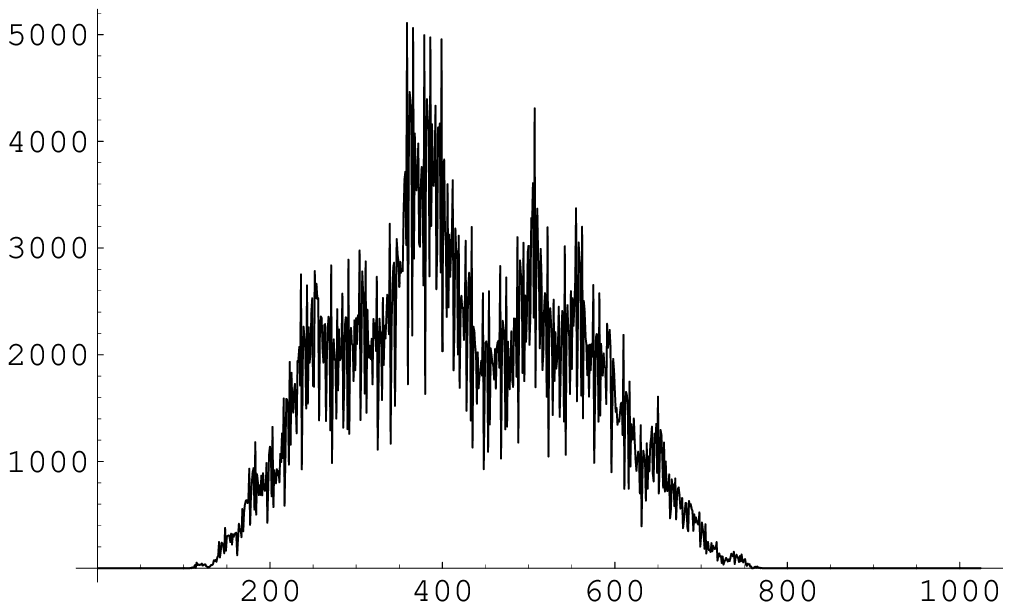}  &
        \epsfxsize=1.7in
        \epsffile{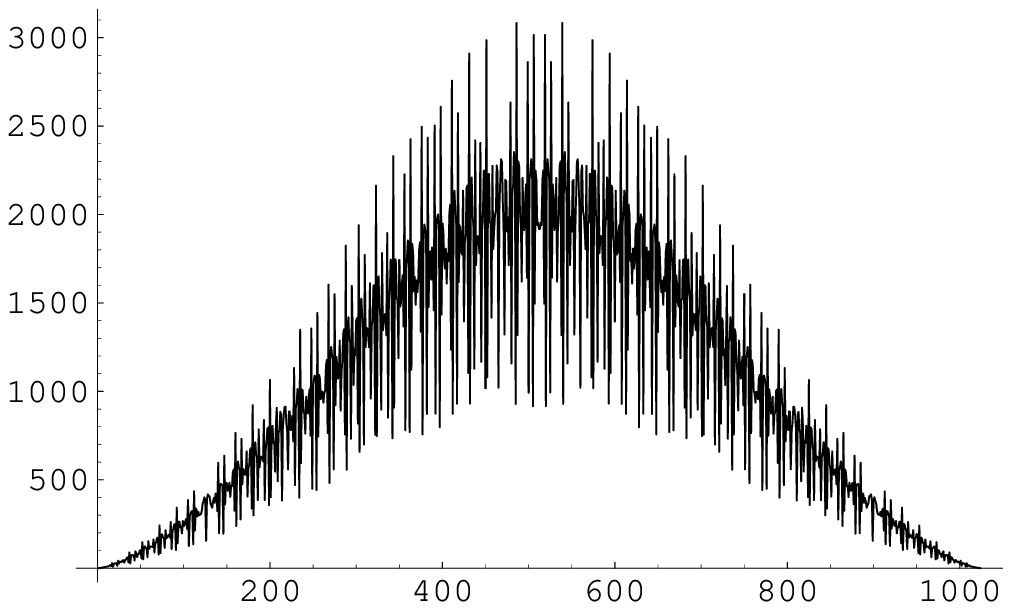} \\ &
\lam = .754877\ldots &
\\ [.3cm] \mbox{simulation (a)} & \mbox{simulation (b)}
  & \mbox{Bernoulli convolution}
\end{array}$
\end{center}
\end{figure}


\section{Properties of the measure}

Here we prove Theorem~\ref{thm-fourier} and the corollaries.
Recall that
$T=\{1,\ldots,\ell\}^*$
is the regular $\ell$-ary tree, with $\ell\ge 2$. 
Let $\eta= \sum_{d\in D} p_d \delta_d$ be a probability distribution on a
finite set $D$ and suppose that $\{a_v\}_{v \in T}$ are i.i.d.\ with the
distribution of $\eta$.
Fix $\lam\in (0,1)$ and
let $\mu$ be the random measure on $\R$ arising from the BRW with
steps of size $\lam^n$ at time $n$ and the behavior of particles governed by
$\eta$.

We can also define $\mu$ using approximating measures at level $n$.
Let $\Ak = \{1,\ldots,\ell\}$.
Define a random measure $\mu_n$ by
\begin{equation} \label{mun}
\mu_n := \sum_{v \in \Ak^n} \frac{1}{\ell^n} \delta_{f(v)}
\ \ \mbox{where}\ \ f(v) := \sum_{j=1}^{n} a_{v|j} \lam^{j}.
\end{equation}

\begin{theorem}[Arbeiter \cite{Arb2}] \label{th-ar1}
The measures $\mu_n$ converge weakly to a random
probability measure $\mu$ almost surely.
\end{theorem}

Now consider the space Sim($\R$) of contracting similitudes on $\R$.
Let $\phi(z) = \lam x+z$ be a map from $\R$ to Sim($\R$). Now define a
probability measure $\Phi:= \otimes_{i=1}^\ell (\eta_i \circ \phi^{-1})$ on the
Borel $\sigma$-algebra of Sim($\R$)$^\ell$, where $\eta_i$ are independent copies
of $\eta$. The measure $\Phi$ is the distribution of a random vector
$(\Fk_1,\ldots,\Fk_\ell)$ of similitudes on $\R$ with contraction ratio
$\lam$ and translation vectors distributed according to $\eta$.

\begin{theorem}[Arbeiter \cite{Arb2}] \label{th-ar2} 
The random measure $\mu$ is stochastically self-similar with respect to $\Phi$.
More precisely,
\begin{equation} \label{eq-ss}
\mu \stackrel{d}{=} \sum_{i=1}^\ell \frac{1}{\ell} \left( \mu^{(i)} \circ
\Fk_i^{-1}\right),
\end{equation}
where $(\Fk_1,\ldots,\Fk_\ell)$ is a random vector of similitudes distributed
according to $\Phi$, and $\mu^{(i)},\ i=1,\ldots,\ell$ are i.i.d.\ copies of
$\mu$ independent of $(\Fk_1,\ldots,\Fk_\ell)$.
\end{theorem}

Next we derive the ``pure types law'' for our measure, which is analogous to the
classical deterministic case.

\begin{proposition} \label{prop-pure} The measure $\mu$ is either a.s.\
absolutely continuous, or a.s.\ pure singular with respect to 
the Lebesgue measure.
\end{proposition}

\begin{proof}
Iterating (\ref{eq-ss}) we obtain for any $n\ge 1$:
$$
\mu \stackrel{d}{=} \sum_{|v|=n} \frac{1}{\ell^n}
\left( \mu^{(v)} \circ
\Fk_v^{-1}\right),
$$
where $\Fk_v$ are some random non-degenerate similitudes $\Fk_v$
(all having the same contraction rate $\lam^n$) and $\mu^{(v)},\ |v|=n$,
are i.i.d.\ copies of $\mu$ independent of what happens in our tree
in the levels $|v|<n$. Consider $\mu_{\rm sing}$, the singular part of $\mu$.
Taking the singular part of the last equation yields
that $\mu_{\rm sing}(\R)$ is independent of the outcome of any finite
number of the random variables on the tree. Thus, $\{\mu_{\rm sing}(\R)>c\}$
is a 0-1 event for any $c>0$, hence there exists $c$ such that
$\mu_{\rm sing}(\R)=c$ almost surely. If $c=0$, then $\mu$ is a.s.\
absolutely continuous. Otherwise, $(1/c) \mu_{\rm sing}$
is a stochastically self-similar probability measure satisfying the equation
(\ref{eq-ss}). By the Uniqueness Theorem (see \cite{Arb1} or 
\cite[Th.\,3.1]{HR}), it equals
$\mu$ in distribution, so $\mu$ is a.s.\ singular. \end{proof}

\smallskip

\begin{proof}[Proof of Theorem~\ref{thm-fourier}]
We follow Bluhm
\cite{Bluhm} closely, up to a point. Fix $n\in \Nat$, then
by (\ref{mun}),
\begin{eqnarray*} \E|\widehat{\mu_n}(t)|^2 & = &
\E\Bigl| \sum_{v\in \Ak^n} \ell^{-n} \exp(itf(v))\Bigr|^2 \\
& = & \sum_{v,w\in \Ak^n} \ell^{-2n} 
\E \exp(it(f(v)-f(w))).\end{eqnarray*}
Observe that
$$
f(v)-f(w) = \sum_{j=k+1}^n \lam^{j}(a_{v|j} - a_{w|j}),
$$
where $k= |v\wedge w|$ and $v\wedge w$ is the longest common initial
segment of $v$ and $w$. Since $a_{v|j}, a_{w|j}$ are
i.i.d.\ with the distribution of $\eta$ for $j=k+1,\ldots,n$, we obtain
\begin{eqnarray*} \E \exp(it(f(v)-f(w))) & = &
\prod_{j=k+1}^n \E \exp(it\lam^{j} a_{v|j})\,
\E \exp(-it\lam^{j} a_{w|j}) \\
& = & \prod_{j=k+1}^n\widehat{ a_{v|j}}(t\lam^{j})
\widehat{a_{w|j}}(-t\lam^{j}) \\
& = & \prod_{j=k+1}^n |\widehat{\eta}(t\lam^{j})|^2,
\end{eqnarray*}
by the definition of the Fourier transform.
Therefore,
\begin{eqnarray}
\E |\widehat{\mu_n}(t)|^2 & = &
\ell^{-2n} \sum_{v,w\in \Ak^n} \prod_{j=|v\wedge w|+1}^n
|\widehat{\eta}(t\lam^{j})|^2 \nonumber \\
& = & \ell^{-2n} \sum_{v\in \Ak^n}\ \sum_{k=0}^n \
\sum_{w:\,|v\wedge w|=k} \ \prod_{j=k+1}^n
|\widehat{\eta}(t\lam^{j})|^2\nonumber \\
& = & \ell^{-n} \Bigl(\sum_{k=0}^{n-1} (\ell-1)\ell^{n-k-1} \prod_{j=k+1}^n
|\widehat{\eta}(t\lam^{j})|^2 + 1 \Bigr).
\label{eq1}
\end{eqnarray}
In the last line we used that $\#\{w\in \Ak^n:\,|v\wedge w|=k\}=
(\ell-1)\ell^{n-k-1}$
for $k=0,\ldots,n-1$, and $|v\wedge w|=n$ only for $v=w$.
Since $|\widehat{\eta}(t\lam^{j})|\le 1$, 
we obtain from (\ref{eq1})
and (\ref{eq-fourier}):
\begin{eqnarray*}
\E |\widehat{\mu_n}(t)|^2 & \ge & \ell^{-n}
\Bigl( (\ell-1) \sum_{k=0}^{n-1} \ell^{n-k-1}
\prod_{j=1}^n |\widehat{\eta}(t\lam^{j})|^2 + 1 \Bigr) \\
& \ge & \prod_{j=1}^n |\widehat{\eta}(t\lam^{j})|^2\ge |\widehat{\nu}(t)|^2.
\end{eqnarray*}
\begin{sloppypar}
By Theorem~\ref{th-ar1}, for every $t$ we have
$|\widehat{\mu}(t)|^2 = \lim_{n\to\infty}
|\widehat{\mu_n}(t)|^2$ hence $\E |\widehat{\mu}(t)|^2 =
\lim_{n\to\infty}\E |\widehat{\mu_n}(t)|^2$ since the Fourier transforms are
bounded by 1 and the expectation is over a finite measure. This proves
part (i) of the theorem.
\end{sloppypar}

(ii) Denote $\theta:= \lam^{-1}$. Let $s$ be a natural number,
$1\le s\le n-1$, suppose that $t\in
[\Th^{s-1},\Th^s]$, and estimate (\ref{eq1}) from above:
\begin{eqnarray*}
& \E |\widehat{\mu_n}(t)|^2 & \le \\ & \le &
\ell^{-n} \Bigl(\sum_{k=0}^{s-1} (\ell-1)\ell^{n-k-1} \prod_{j=k+1}^n
|\widehat{\eta}(t\lam^{j})|^2 + (\ell-1)\sum_{k=s}^{n-1}\ell^{n-k-1}+
1 \Bigr) \\ & = & (\ell-1) \sum_{k=0}^{s-1} \ell^{-k-1} \prod_{j=k+1}^n
|\widehat{\eta}(t\lam^{j})|^2 + \ell^{-s}.
\end{eqnarray*}
Letting $n\to\infty$ and using (\ref{eq-fourier}) we obtain
\begin{equation} \label{eq2}
\E |\widehat{\mu}(t)|^2 \le \ell^{-s} + (\ell-1) \sum_{k=0}^{s-1} \ell^{-k-1}
|\widehat{\nu}(t\lam^k)|^2.
\end{equation}
Next we multiply by $t^{2\gam}$ and
integrate over $[\Th^{s-1},\Th^s]$, keeping in mind that
$\lam = \Th^{-1}$ and
$$
\int_{\Th^{s-1}}^{\Th^s} |\widehat{\nu}(t\lam^k)|^2 t^{2\gam}\,dt =
\Th^{k(1+2\gam)}
\int_{\Th^{s-k-1}}^{\Th^{s-k}} |\widehat{\nu}(t)|^2 t^{2\gam}\,dt
$$
to obtain
\begin{eqnarray*}
\int_{\Th^{s-1}}^{\Th^s} \E |\widehat{\mu}(t)|^2 t^{2\gam}\,dt & \le &
\frac{\Th^{s(1+2\gam)}-\Th^{(s-1)(1+2\gam)}}{\ell^s(1+2\gam)} \\[1.1ex]
& + & (\ell-1) \sum_{k=0}^{s-1}
\frac{\Th^{k(1+2\gam)}}{\ell^{k+1}}\int_{\Th^{s-k-1}}^{\Th^{s-k}}
|\widehat{\nu}(t)|^2 t^{2\gam}\,dt.
\end{eqnarray*}
Summing over $s=1,2,\ldots,$ using that $\Th^{1+2\gam}< \ell$, and exchanging
the order of summation yields for some $C_1,C_2$:
$$
\int_1^\infty \E |\widehat{\mu}(t)|^2 t^{2\gam}\,dt \le C_1 + C_2 \int_1^\infty
|\widehat{\nu}(t)|^2 t^{2\gam}\,dt,
$$
and this implies (\ref{ineq2}). \end{proof}

\smallskip

\begin{proof}[Proof of Corollary~\ref{cor2}] (i) If $\mu$ is not almost
surely singular,
then it is almost surely
absolutely continuous by Proposition~\ref{prop-pure}, hence
by the Riemann-Lebesgue Lemma we have $\lim_{|t|\to\infty} |\muhat(t)|^2=0$
almost surely. Since $|\muhat(t)|\le 1$ for all $t$, this being a probability
measure, by the Lebesgue Dominated Convergence Theorem we have
$\lim_{|t|\to\infty} \E|\muhat(t)|^2=0$, contradicting the assumption on $\nu$
in view of Theorem~\ref{thm-fourier}(i).

\begin{sloppypar}
(ii) By Plancherel's Theorem, $\nu$ has a density in $L^2$ if and only if
$\int_{\R} |\nuhat(t)|^2\,dt < \infty$. Now the claim follows from
Theorem~\ref{thm-fourier}(i).
\end{sloppypar}

(iii) and (iv) are  
immediate from Theorem~\ref{thm-fourier}(ii) and Fubini Theorem.
\end{proof}

\smallskip

\begin{proof}[Proof of Corollary~\ref{cor-meas}]
(i) It is well-known that $\lim_{|t|\to\infty} |\widehat{\nu}(t)|\ne 0$
when $1/\lam>m$ is Pisot. This is due to Erd\H{o}s \cite{Erd1} for $m=2$,
and the proof easily extends to arbitrary $m$ (see \cite{BG}).
Now the claim follows from Corollary~\ref{cor2}(i).

(ii) By the result of Simon and T\'oth \cite{ST}, which extended \cite{Sol}
to the case $m>2$, the self-similar measure $\nu$ is a.c.\ with a density in
$L^2$ for a.e.\ $\lam\in (\frac{1}{m},1)$.
Now the claim follows from Corollary~\ref{cor2}(iii).
\end{proof}

\smallskip

\begin{proof}[Proof of Corollary~\ref{cor3}]
(i) follows from \cite{Garsia} and Corollary~\ref{cor2}(iii).

(ii) follows from \cite[Lem.\,5]{PS} (see also \cite{Erd2} for the classical,
but less sharp result) and  Corollary~\ref{cor2}(iv).

(iii) follows from \cite[Section 5]{PS}. \end{proof}


\section{Intervals in the support}

In this section we prove Theorem~\ref{th-int}. 
Here we return to our most basic set-up: BRW with two equally likely digits
on the binary tree. As already mentioned, it is convenient to use
the digits 0,1, so we have i.i.d.\ random variables
distributed as $\half \delta_0  + \half \delta_1$, at each vertex of
the tree.

First we introduce some notation. For $a = a_1 a_2\ldots$, a finite
or infinite word in the alphabet $\{0,1\}$, denote
$$
\xi(a) = \xi(a_1 a_2\ldots) = \sum_{j=1}^\infty a_j \lam^{j} 
\ela .a_1a_2\ldots
$$
Recall that $g=\frac{\sqrt{5}-1}{2}$, so that $1=g + g^2 = \sum_{j=2}^\infty
g^j$. Since $\lam>g$, there exists $\ell\ge 3$ such that
\begin{equation} \label{eq-star}
1 < \lam^2 + \cdots + \lam^\ell.
\end{equation}
We fix $\ell$ for the rest of the proof.
Let
$$
\Uk:= (\alpha,\beta),\ \ \mbox{where}\ \
\alpha = \frac{\lam^\ell}{1-\lam^\ell}\,,\
\ \ \ \beta = \frac{\lam}{1-\lam} - \frac{\lam^\ell}{1-\lam^\ell}\,.
$$
Note that
$$
\alpha \ela .(0^{\ell-1}1)^\infty,\ \ \ \beta \ela .(1^{\ell-1}0)^\infty.
$$
For $a\in \{0,1\}^n$, let
$$
\Uk_a := \xi(a) + \lam^n \Uk.
$$

\begin{lemma} \label{lem-greedy} There exists $c>0$ such that for any
subinterval $J\subset \Uk$, with $|J| \le c$, there exist $a,a' \in
\{0,1\}^{\ell-1}$ satisfying
$$
J \subset \Uk_{a0} \cap \Uk_{a'1}.
$$
\end{lemma}

\begin{proof}
The idea comes from \cite[Th.\,3]{EJK} which shows
that given any $x\in \Uk$, we can fix an arbitrary sequence $\{a_{kl}\}_{k\ge 1}
\in \{0,1\}^\infty$ and obtain an expansion $x \ela .a_1 a_2\ldots$ with
the digits 0,1 by the ``greedy algorithm.'' It is enough to show that
\begin{equation} \label{claim1}
\Uk \subset \Bigl( \bigcup_{a \in \{0,1\}^{\ell-1}} \Uk_{a0} \Bigr) \cap
\Bigl( \bigcup_{a \in \{0,1\}^{\ell-1}} \Uk_{a1} \Bigr).
\end{equation}
Indeed, since the intervals are open, it will follow that there exists a
positive $c$ as desired.

Let $a\in \{0,1\}^{\ell-1}$ be non-maximal, that is, $a \ne 1^{\ell-1}$.
We claim that there exists $a'\in \{0,1\}^{\ell-1}$ such that
$\xi(a) < \xi(a')$ and
\begin{equation} \label{claim2}
\Uk_{aj} \cap \Uk_{a'j} \ne \es,\ \ \mbox{for}\ j=0,1.
\end{equation}
This is equivalent to showing that
\begin{equation} \label{claim3}
\xi(a')-\xi(a) < |\Uk_{aj}| = \lam^{-\ell} |\Uk|.
\end{equation}
If $a$ ends with $0$, consider $a'$ which ends with $1$, but 
otherwise agrees with $a$. Then $\xi(a')-\xi(a) = \lam^{\ell-1}$.
If $a$ ends with $1$, then $a$ ends
with $01^p$ for some $1\le p\le \ell-2$ since $a$ is non-maximal.
Consider the greedy expansion of $1$ in base $\lam$:
$$
1=d_1\lam + d_2 \lam^2 + \cdots
$$
It has the property that $d_j\in \{0,1\}$ for all $j$ and 
$1-\sum_{j=1}^n d_j \lam^j < \lam^n$ for all $n\ge 1$. Let
$$
w = 1,(1-d_1),\ldots,(1-d_p)
$$
and consider $a'\in \{0,1\}^{\ell-1}$ which ends with $w$,
but otherwise agrees with $a$. Then
\begin{eqnarray*}
\xi(a')-\xi(a) & = & \lam^{(\ell-1)-(p+1)} (\xi(w)-\xi(1^p))\\
& = & \lam^{\ell-p-1}
\Bigl(1 - \sum_{j=1}^p d_j \lam^j\Bigr) \in (0, \lam^{\ell-1}).
\end{eqnarray*}
In both cases we obtain that $\xi(a')-\xi(a) \le \lam^{\ell-1}$. 
Note that $1 - \lam^{\ell} < \lam^2 + \cdots + \lam^{\ell} - \lam^{\ell+1}$
by (\ref{eq-star}), hence
$$
\lam^{-1} < \frac{\lam+ \cdots + \lam^{\ell-1} - \lam^\ell}{1-\lam^\ell} =
 \frac{\lam}{1-\lam} - \frac{2\lam^\ell}{1-\lam^\ell} = |\Uk| = \lam^{-\ell}
|\Uk_{aj}|.
$$
Therefore, 
$$
\xi(a')-\xi(a) \le \lam^{\ell-1} < |\Uk_{aj}|,
$$
which verifies (\ref{claim3}) and hence (\ref{claim2}). It remains to
note that $\alpha$ (the left endpoint of $\Uk$) is the left endpoint
of $\Uk_{0^{\ell-1}1}$ and is covered by $\Uk_{0^\ell}$; similarly,
$\beta$ is the right endpoint of $\Uk_{1^{\ell-1}0}$ and is covered by
$\Uk_{1^\ell}$. This, together with (\ref{claim2}), implies (\ref{claim1}), and
the lemma is proved. 
\end{proof}

\begin{corollary} \label{cor-greedy} There exists a constant
\begin{equation} \label{c-cond}
0 < c < \lam^\ell|\Uk|/4
\end{equation}
such that for any $n\ge 1$ and any $a \in \{0,1\}^{(n-1)\ell}$,
if $J\subset \Uk_a$, with $|J| \le c\lam^{n\ell}$, then there exist $a',a'' \in
\{0,1\}^{\ell-1}$ satisfying
$$
J \subset \Uk_{aa'0} \cap \Uk_{aa''1}.
$$
\end{corollary}

\begin{proof}
For $n=1$ this is just Lemma~\ref{lem-greedy} (of course,
we can always impose an upper bound on $c$). For $n>1$ this follows from
Lemma~\ref{lem-greedy} and the definitions by rescaling,
with the same $c$. \end{proof}

\smallskip

\begin{proof}[Proof of Theorem~\ref{th-int}]
We are going to use a variant of the argument from \cite{BK}. In order to set
it up, we need to consider another family of intervals, which we denote $J_w$
and which should not be confused with the intervals $\Uk_a$.
Fix $p\in \Nat$ such that
\begin{equation} \label{p-cond}
p^{-1} < \min\{c(1-\lam)/\lam, \lam^\ell\}
\end{equation}
where $c$ is from Lemma~\ref{lem-greedy}.
We subdivide the interval $[0,\frac{\lam}{1-\lam}]$ into $p$ closed
subintervals of equal length and denote them $J_i,\ i=0,\ldots,p-1$.
Then continue the subdivision and denote the intervals of level $n$
by $J_w$ for $ w\in \{0,\ldots,p-1\}^n$. We have
$ |J_w| = \frac{\lam}{1-\lam} p^{-|w|} $, 
so in view of (\ref{p-cond}) and (\ref{c-cond}),
\begin{equation}\label{J-cond}
|J_w| < c\lam^{\ell(n-1)} < |\Uk_a|/4 \ \ \ \mbox{for}\ |w|=n \ \mbox{and}\
|a|=n\ell.
\end{equation}
For a vertex $\sig$ of the binary tree denote by $a(\sig)$
the (random) finite word of 0's and
1's which we see on the path from the root to $\sig$. 
Fix $\eps$ such that
\begin{equation} \label{eps-cond}
0 < \eps < 2\lam-1\ \ \mbox{and}\ \ (1+\eps)^\ell < 2.
\end{equation}
Now we follow the argument of \cite[pp.\,1046-47]{BK} very closely.
For each $n\ge 1$ and $w\in \{0,\ldots,p-1\}^n$ define the event
\begin{eqnarray}
H_n(w)
 & := & \{\mbox{there are at least}\ (1+ \eps)^{n\ell}\ \mbox{vertices}
\nonumber \\ \label{def-event}
       &  & \sigma,\ \mbox{with}\ |\sig|=n\ell,\
            \mbox{such that}\ J_w \subset \Uk_{a(\sig)}\}.
\end{eqnarray}
We are going to show that
\begin{equation} \label{eq-sum}
\sum_{n=1}^\infty \
\sum_{w_1\ldots w_{n+1}}
\Prob\{H_n(w_1\ldots w_{n})
\setminus H_{n+1}(w_1\ldots w_{n+1})\} < \infty
\end{equation}
where the second sum is over all sequences $w\in \{0,\ldots,p-1\}^{n+1}$.
First we explain how this implies the desired result.
By Borel-Cantelli, (\ref{eq-sum}) implies that almost surely for all
$n$ sufficiently large ($n\ge N$ where $N$ is random),
$H_n(w_1\ldots w_{n})
\setminus H_{n+1}(w_1\ldots w_{n+1})$ does not occur.
Let us fix a random configuration (the choice of 0's and 1's), so that
$N$ is now fixed. We write $\sim$ to indicate equality up to a 
(multiplicative) positive
constant independent of $n$.
We have $2^{n\ell}$
intervals $\Uk_{a(\sig)}$ of level $n\ell$,
each of length $\sim\lam^{n\ell}$, whose union is contained in
$[0,\frac{\lam}{1-\lam}]$. Note that many of them will likely
coincide; they are counted with multiplicity.
By the pigeon-hole principle,  there are at least $\sim 2^{n\ell}
\lam^{n\ell}$ intervals $\Uk_{a(\sig)}$ of level $n\ell$ with a common
intersection longer than $|\Uk_{a(\sig)}|/2$. For large $n$ the number of these
intervals exceeds
$(1+\eps)^n$ since $2\lam > 1+\eps$, hence there is an interval
$J_w$, with $|w|=n$, contained in their intersection by (\ref{J-cond}).
Thus, the event $H_n(w)$ occurs, and we can assume that $n\ge N$. By the
choice of $N$, the events $H_j(ww')$, with $j=|ww'|$,
occur for all finite extensions
of the word $w$. This implies that a.s.
$$
J_w \subset \bigcap_{n=1}^\infty
\bigcup_{j=n}^\infty \bigcup_{|\sig|=j} \Uk_{a(\sig)}.
$$
It is easy to see that the right-hand side of the last formula is contained
in $S$: it consists of points which can be approximated by
$\xi(a(\sig))$ for $|\sig|$ arbitrarily large, and $S$ is compact.
Thus, $S$ contains an interval a.s. The same argument,
of course, implies that a.s.\ there is an interval in the ``cylinder''
of $S$ obtained by taking a subtree from any given vertex, and since such
cylinders are dense in $S$, it follows that $S$ is the closure of its
interior.

It remains to verify (\ref{eq-sum}). We continue to follow the scheme of
\cite[Section 5]{BK}.
Assume that $H_n(w)$ occurred, so that there exist at least
$(1+\eps)^{n\ell}$ vertices $\sig_1,\ldots,\sig_r$ at the level
$n\ell$ such that $J_w\subset \Uk_{a(\sig_j)}$, $1 \le j \le r$. Let
$w'=ww_{n+1}$. We have $J_{w'} \subset J_w$ and $|J_{w'}| < c\lam^{n\ell}$.
We can apply Corollary~\ref{cor-greedy} for $J=J_{w'}$ and each
$\Uk_{a(\sig_j)}$ to conclude that there exist
$a'_j,a''_j \in \{0,1\}^{\ell-1}$
satisfying
$$
J_{w'} \subset \Uk_{a(\sig_j)a'_j0} \cap \Uk_{a(\sig_j)a''_j1}.
$$
Conditionally
on $H_n(w)$, the following random variables are independent of each other
for $j=1,\ldots,r$:
\begin{eqnarray*}
U_n(\sig_j) & := & \{\mbox{number of descendants $\tau$ of $\sig_j$ in the}\
                  \mbox{$(n+1)\ell$-th level} \\
 & &             \mbox{with the property that the word seen on the path from}
         \\ & & \mbox{$\sig_j$ to $\tau$ is either $a'_j0$ or $a''_j1$} \}.
\end{eqnarray*}
Moreover,
$$
U_n(\sig_j) \ge 0\ \ \mbox{and}\ \ \E\{U_n(\sig_j)\,|\, H_n(w), \sig_1,\ldots,
\sig_r\} = 2.
$$
Indeed, the expected number of times  to read off
a specific word $W$ in the tree of
descendants of any vertex from level $n\ell$ down to level $(n+1)\ell$
equals the number of paths times
the probability of seeing $W$ on a given path. This gives $2^\ell\cdot 
2^{-\ell}=1$,
(recall that we are in the unbiased case), but since two words are good for us,
we get 2.  Note that
$U_n(\sig_j)\le 2^\ell$, so these random variables are
uniformly bounded, and one easily checks
that the variance satisfies
$$
0 < C_1 \le {\rm Var}(U_n(\sig_j)) \le C_2 < \infty
$$
for some $C_1, C_2$ independent of $n$ and $j$. Therefore, by the
Large Deviation Estimate (see e.g.\ \cite[Section 7.4]{renyi}), in view of 
$(1+\eps)^\ell<2$,
there exists a constant $C_3 = C_3(\eps)$ such that
$$
\Prob\Bigl\{\sum_{j=1}^r U_n(\sig_j) \le r(1+\eps)^\ell\,|\, 
H_n(w),\sig_1,\ldots,\sig_r
\Bigr\} \le 2\exp(-C_3 r).
$$
However, if $H_n(w)$ occurs, and hence $r\ge (1+\eps)^{n\ell}$, and
$\sum_{j=1}^r U_n(\sig_j) > r(1+\eps)^\ell \ge (1+\eps)^{(n+1)\ell}$, 
then there are
at least $(1+\eps)^{(n+1)\ell}$ vertices $\tau$ in the $(n+1)\ell$-th level
such that $J_{w'} \subset \Uk_{a(\tau)}$, that is, $H_{n+1}(w')$ occurs.
This implies
$$
\Prob\{H_{n+1}(w')\ \mbox{fails}\,|\,H_n(w)\} \le 2\exp(-C_3(1+\eps)^{n\ell}),
$$
and (\ref{eq-sum}) follows. \end{proof}

We do not know whether the set $S$ contains an interval a.s.\
for $\lam\in (\half,g)$. However, if we can
find an interval which is covered with multiplicity at least two by its images
on certain level, then the proof above goes through without any changes.
We state this precisely for future reference.

\begin{proposition} \label{prop-new}
Let $\lam\in (\half,1)$ and suppose that there exists a  nonempty open
interval $\Uk$ and $\ell\in \Nat$ satisfying the condition
(\ref{claim1}).
Let $\mu$ be the BRW with steps
of size $\lam^{n}$, with the i.i.d.\ random variables on the binary tree
distributed as $\half(\delta_0 + \delta_1)$.
Then $S=\supp(\mu)$
has nonempty interior, and is the closure of its interior, a.s.
\end{proposition}

\smallskip

\begin{proof}[Proof of Proposition~\ref{prop-discon}]
As explained in the Introduction, if $\Ek_\lam(x)$ has zero logarithmic
capacity, then $x\not\in S$ a.s. We can choose a countable dense set in $I$
with this property, and the claim follows. \end{proof}

\smallskip

\begin{proof}[Proof of Corollary~\ref{cor-gold}]
It is known (easy to see, folklore, see e.g.\ \cite[App.\,A]{SiVer})
that for $\lam=g$ every point $x\equiv ng^{-1}$ (mod 1) 
has countably many expansions
in base $g$. Since this is a countable dense set, Proposition~\ref{prop-discon}
applies. \end{proof}

Now we turn to the proof of Proposition~\ref{prop-easy}.
Let $I =
[0,\frac{\lam}{1-\lam}]$ and recall that $S\subset I$.
For $v\in T$, let $I(v) =
f(v) + \lam^{|v|} I$. We have
\begin{equation} \label{eeq1}
S = \bigcap_{n=1}^\infty S(n),\ \ \ \mbox{where}\ \ S(n):=\bigcup_{|v|=n} I(v).
\end{equation}
The sets $S(n)$ form a decreasing nested family.

\smallskip

\begin{proof}[Proof of Proposition~\ref{prop-easy}]
Let $K_n$ be the number of distinct words of length $n$ which we see from
the root, that is, $K_n = \#\{a(v):\ |v|=n\}$. The probability of seeing any 
given word is $\sim n^{-1}$, since this is the probability of survival of
a critical branching process (see e.g.\ \cite[Th.\,I.9.1]{AN}).
Summing over all possible words
we obtain $\E(K_n) \sim n^{-1}2^n$.
It follows that 
$\E|S(n)| \le C\lam^n \E(K_n) \to 0$, as $n\to\infty$, for $\lam\le \half$.
Clearly, $\E|S| = \lim_{n\to\infty}  \E|S(n)|=0$, hence $|S|=0$ almost surely.
\end{proof}


\section{The support is fractured}

Here we prove Theorem~\ref{th-fract}.
We will use the following fact about Pisot numbers.

\begin{lemma}[Garsia \cite{Garsia}] \label{lem-gar}
For every $\lambda=1/\theta$, where $\theta$ is Pisot,
there is a constant $0<c_1<1$ such that if $a_1,..,a_n$ and $b_1,..,b_n$ are
sequences of $0,1$ and
$\sum_{i=1}^n a_i \lambda^i \ne \sum_{i=1}^n b_i \lambda^i$, then
$|\sum_{i=1}^n a_i \lambda^i - \sum_{i=1}^n b_i \lambda^i| \ge c_1 \lambda^n$.
\end{lemma}

\begin{proof}[Proof of Theorem~\ref{th-fract}]
We denote by $\ba\in \{0,1\}^T$ the outcome of the lotteries on the
tree, that is, $\ba = \{a_v:\ v\in T\}$.
We also write $T_n$ for the set of vertices at level $n$ and use similar 
notation for subgraphs of $T$.

Let $q$ be a rational number, and let $\Pk_q$ be the event that
$q\not\in S$. It is enough to prove that the endpoints of the
component of $\R\setminus S$ containing $q$ have the desired
property a.s., conditioned on $\Pk_q$. We will work with the left
endpoints, since the right endpoints are treated exactly the same
way. So fix $q$, assume $\Pk_q$ holds, and let $(\alpha,\beta)$ be
the component of $\R\setminus S$ containing $q$.

Let  $n_0$ be the smallest integer such that $q\not\in S(n_0)$. By
(\ref{eeq1}), there exists such an integer. Note that conditioned on
$\Pk_q$ and on the value of $n_0$, the distribution of $\ova$ for
levels greater then $n_0$ is still the same product measure. This is
because the event $q\not\in S(n_0)$ depends only on the first $n_0$
levels of $\ova$, and $\Pk_q$ is simply the union of the corresponding 
cylinder sets.

For any $n\ge n_0$, let
$$ L^{(n)}=\{ v\in T_n : f(v)<q\} ,$$
and let $m^{(n)}=\max \{f(v) : v\in L^{(n)}\}$. Let
$$ \Mk^{(n)}=\{v\in L^{(n)} : f(v)=m^{(n)}\}$$
be the set of vertices where this maximum is achieved. Members of
$\Mk^{(n)}$ are called the {\em maximal vertices of level $n$}. Of
course, $\Mk^{(n)}$ is a random subset of $T_n$.

For $v\in T$ consider the subtree
$$ T_v^+=T_v^+(\ova) = \{vw:\ w\in \{1,2\}^*,\ a_{vw_1} = a_{vw_1 w_2} =
\cdots =
a_{vw}=1\}.
$$
In other words, $T_v^+$ is the cluster of 1's from the vertex $v$ down (we
do not make any assumptions on $a_v$). The number of vertices of
$T_v^+$ at level $n$ is a critical Galton-Watson process since the
number of offsprings of a vertex labeled by 1 has the distribution
$\frac{1}{4}\delta_0 + \half \delta_1 + \frac{1}{4}\delta_2$, with
expectation 1. Such a process dies out with probability 1, so $T_v^+$
is a.s.\ finite, for all $v\in T$.

Fix $n\ge n_0$, consider the graph
$$ \Gam^{(n)}= \bigcup\{T_v^+:\ v\in \Mk^{(n)}\}, $$ and let
$$ \tau^{(n)} = \sup\{j:\ \Gam^{(n)}_j \ne \es\}.  $$
We have $\tau^{(n)}< \infty$ a.s.
Fix $\ell\in \Nat$ such that
\begin{equation}\label{def-ell} \frac{\lam^{\ell+1}}{1-\lam} < c_1,
\end{equation}
where $c_1$ is from Lemma~\ref{lem-gar}.
Consider the event
$$ \Ek^{(n)}:=\left\{\tau^{(n)} < \infty\ \&\
|\Gam^{(n)}_{\tau^{(n)}-j}| = 2^{\ell-j}, \ j=0,\ldots,\ell\right\}.  $$
In other words, $\Ek^{(n)}$ occurs whenever
$\Gam^{(n)}_{\tau^{(n)}-\ell} = \{\sig\}$ for some vertex $\sig$ and
all its $2^{\ell+1} - 2$ descendants down to the level $\tau^{(n)}$ are
labeled by 1's, after which $\Gam^{(n)}$ dies out.

\begin{lemma} \label{lem-est}
$
\exists\,c_2>0,\ \forall\,n\in \Nat,\ \Prob(\Ek^{(n)}) \ge c_2.
$
\end{lemma}

\begin{proof} This is written in measure-theoretic language.
Denote $\Omega = \{0,1\}^T$ and recall that $\Prob = (\half,\half)^T$,
so that $(\Omega,\Prob)$ is the probability space for our lotteries on the tree.
We show the dependence on $\ova$ in our notation.
Let $\Omega_0 = \{\ova\in \Omega:\ \tau^{(n)}(\ova) <
\infty\}$. The fact that all $T_v^+$, for $v\in \Mk^{(n)}$,
die out a.s.\ means that
$\Prob(\Om_0) = 1$.

Consider the following transformation on $\Om_0$. Let $\ova \in
\Om_0$ and choose $w$ to be the rightmost (i.e.\ greatest in the
lexicographic order) vertex of $\Gam^{(n)}$ in the level 
$\tau^{(n)}(\ova)$, that is, just before it dying out. By the
definition of $\tau^{(n)}(\ova)$ we have $a_{w 1} = a_{w 2} = 0$.
Then $\ova'\in \Om_0$ is defined as follows:
\begin{eqnarray*}
a'_{w1 u} = 1 & \forall\,u\in \{1,2\}^*, & 0\le |u|\le \ell; \\
a'_{w 1 u} = 0 & \forall\,u\in \{1,2\}^{\ell+1}. &
\end{eqnarray*}
For all other vertices $v$ we let $a'_v = a_v$.

Observe that $\tau^{(n)}(\ova') = \tau^{(n)}(\ova) + \ell + 1$,
and $\ova'\in \Ek^{(n)}$ by construction (note that $\sig = w1$ is the
vertex indicated after the definition of $\Ek^{(n)}$).
Let $\Om' = \{\ova':\ \ova \in \Om_0\}$. It is enough to prove that
\begin{equation} \label{eqe2}
\Prob(\Om') \ge c_2 > 0
\end{equation}
for some constant $c_2$ which does not depend on $n$. This
follows from the definition of $\Prob$ as the product measure on $\Om$
and the fact that given $\ova'$ we can recover $\ova$ except for the
descendants of $\sig 1$ in the levels $\tau^{(n)}(\ova)+1,\ldots,
\tau^{(n)}(\ova)+\ell+2$. Thus, we can take $c_2 = 2^{-2^{\ell+2}}$.
\end{proof}

\smallskip

\begin{lemma} \label{lem-gap}
If $\Ek^{(n)}$ occurs, then $(\alpha - \frac{\lam^n}{1-\lam},\alpha)
\not\subset S$.
\end{lemma}

\begin{proof}
To simplify notation, we write $\Gam = \Gam^{(n)}$ and
$\tau=\tau^{(n)}$. Suppose that $\Ek^{(n)}$ occurred, so that
$\Gam_{\tau-\ell} = \{\sig\}$ for some vertex $\sig$, and all its
$2^{\ell+1} - 2$ descendants down to level $\tau$ are labeled by
1's, after which $\Gam$ dies out.

Recall that $\Gam$ is the cluster of 1's starting from the maximal
vertices at level $n$. Thus, for $j\le \tau^{(n)}$ we have
$\Gam_j=\Mk^{(j)}$. It follows that $\sig$ is the only maximal
vertex at the level $\tau-\ell$. For any $u$, with $|u|\ge \ell$,
\begin{equation} \label{ika1}
f(\sig u) \ge f(\sig) + \lam^{\tau-\ell+1} + \cdots + \lam^{\tau} =
f(\sig) + \lam^{\tau-\ell+1} \frac{1-\lam^\ell}{1-\lam}\,.
\end{equation}
Now suppose $v\in T_{\tau-\ell},\ v\ne \sig$. Then $v$ is not
maximal, hence for any $u\in \{1,2\}^*$, if $f(vu) < q$ then
\begin{equation} \label{ika11}
f(vu) \le f(v) + \frac{\lam^{\tau-\ell+1}}{1-\lam} \le f(\sig) - c_1
\lam^{\tau-\ell}  + \frac{\lam^{\tau-\ell+1}}{1-\lam}
\end{equation}
by  Lemma~\ref{lem-gar}. Observe that the right-hand side of
(\ref{ika11}) is less than the right-hand side of (\ref{ika1}) by
(\ref{def-ell}), hence we have a gap.

Since  $\sigma$ is a maximal vertex at the level $\tau-\ell$, we
have
\begin{equation} \label{ika111}
\alpha = \max(S\cap (-\infty,q)) \in \Bigl[ f(\sig) +
\lam^{\tau-\ell+1} \frac{1-\lam^\ell}{1-\lam}, f(\sig) +
\lam^{\tau-\ell+1} \frac{1}{1-\lam} \Bigr].
\end{equation}
In view of (\ref{ika1}), we have
$[A-\eps,A]\not\subset S$ for 
$A=f(\sig) +\lam^{\tau-\ell+1} \frac{1-\lam^\ell}{1-\lam}$ and $\eps>0$.
Now (\ref{ika111}) and $\tau\ge n$ imply the desired claim
$(\alpha - \frac{\lam^n}{1-\lam},\alpha)\not\subset S$. 
\end{proof}

\smallskip

{\bf Conclusion of the proof.} Since the critical branching process
dies out a.s., we can find a sequence $n_k\uparrow \infty$ such that
$n_{k+1}-n_k > \ell$ and with probability greater than $1-\delta$,
for some fixed $\delta$, for all $k$ and all vertices $\sigma
\in T_{n_k}$, the random subgraph $T^+_\sig$ dies out before level
$n_{k+1}$.
 Let $\Gk$ be the event that this happens (so that $\Prob(\Gk)> 1-\delta$).
Conditioned on this event, the events $\Ek^{(n_k)}$, considered
above, are independent. In view of Lemma~\ref{lem-est}, infinitely
many of the events $\Ek^{(n_k)}$ occur a.s., conditioned on $\Gk$.
By Lemma~\ref{lem-gap}, this implies that $(\alpha-\eps,\alpha)
\not\subset S$ for every $\eps>0$. Since $\delta>0$ is arbitrary,
the claim of the theorem follows.
\end{proof}


\section{Further Directions}

Here we discuss some directions for further research and
mention some open questions (in addition to those listed in Section 1).

\smallskip

{\bf 1.}
The signs on the tree are taken with different probabilities, e.g.\
$(p,1-p)$ on the binary tree. 
Biased Bernoulli convolutions exhibit
``multifractal'' behavior, and the thresholds for absolute continuity and
$L^q$ density no longer coincide, see \cite{PeSo2}. We can expect a similar 
phenomenon for our random model. Combined with the results of \cite{PeSo2},
Corollary~\ref{cor2} 
yields sufficient conditions for the a.s.\ existence of $L^2$ density, 
at least for some
values of $p$. Entropy considerations imply that for all $p\ne \half$,
for some $\lam>\half$ the measure $\mu$ is a.s.\ supported on a set of zero 
Lebesgue measure. The methods of Sections 3 and 4 partially extend.
For instance, if $p>\half$, then for $1/\lam$ Pisot every gap of
the support has gaps accumulating to it from the right.

{\bf 2.} The following model was suggested by K\'aroly Simon: fix $p$
and suppose that the random variables on the vertices are independent
and have either the distribution $p\delta_1 + (1-p) \delta_{-1}$ or
$(1-p)\delta_1 + p \delta_{-1}$, depending on whether the last edge goes
left or right. This can be considered as a $p$-perturbation of the
Bernoulli convolution measure, since for $p=0$ we get the Bernoulli 
convolution (there is no randomness). For $p=\half$ we get the model studied
in this paper. What happens for $0< p < \half$, especially as $p\to 0$?


{\bf 3.} Other graphs (1): instead of the regular rooted tree, consider more
general graphs, e.g.\ 1/4 of the grid $\Z^2$. There are $2^n$
paths of length $n$, but there is less randomness, since there are only
$\sim n^2$ edges (put the signs on edges).

{\bf 4.}  Other graphs (2): consider the same problem on a random tree,
e.g.\ on a Galton-Watson tree (choosing a random tree is part of the model).
This way we recover stochastic self-similarity (in some sense).

{\bf 5.} Other graphs (3): consider the same problem on a random or
deterministic graph of  polynomial growth. Instead of the steps $\lam^n$
consider slowly shrinking steps. 
E.g.\ spherically symmetric trees of growth $r^d$, 
or critical GW tree conditioned
to survive and steps at level $n$ equal to  $\pm n^{-1/2}$. 
Here there is no corresponding
Bernoulli convolution.

{\bf 6.} The random variables $a_v$
on the tree are i.i.d.\ with some a.c.\ distribution, e.g.\
Gaussian. (Note that we first choose these steps, fix them, and then consider
the BRW with steps $a_v \lam^{|v|}$, that is, the contraction ratio is 
deterministic.)
On the binary tree, there will be $2^n$ intervals of size $\sim \lam^n$,
so for $\lam < 1/2$ we 
still get singularity. For $\lam> 1/2$ we should have
a.c.\ almost surely; the ``Pisot numbers effect'' will be lost.
This should follow easily by the methods of \cite{PSS,JPS}.
We do not know what to expect 
for connectedness properties of the support.

\smallskip

{\bf Remark added on May 30, 2007.} Russ Lyons indicated to us a simpler
proof of Corollary 1.2(i) and Corollary 1.3(i). In fact, the expectation of
the random measure $\mu$ in Theorem 1.1 
over all the ``lotteries'' on the tree is 
easily seen to be the self-similar measure $\nu$ defined by (1.1). 
Therefore, if $\nu$ is singular, then $\mu$ is a.s.\ singular, by Fubini
Theorem.

\smallskip

{\bf Acknowledgments.} Thanks to Nikita Sidorov for many helpful comments
concerning expansions in non-integer bases, and to Russ Lyons,
Yuval Peres and K\'aroly
Simon for useful discussions.
This work was done during the visit of 
Boris Solomyak to the Weizmann Institute of Science, as 
Rosi and Max Varon Visiting Professor in the Autumn of 2005.
He is grateful to the Department of Mathematics at WIS
for its hospitality. 


\end{document}